# Unity in Major Themes – Convergence vs. Arbitrariness in the Development of Mathematics[1]

## Bernhelm Booß-Bavnbek (Roskilde University, Denmark)[2] and

## Philip J. Davis (Brown University, Providence, Rhode Island, U.S.A.)[3]

*To the memory of* **Gian-Carlo Rota** *(April 27, 1932 – April 18, 1999)*

Abstract. We describe and explain the desire, common among mathematicians, both for unity and independence in its major themes. In the dialogue that follows, we express our spontaneous and considered judgment and reservations; by contrasting the development of mathematics as a goal-driven process as opposed to one that often seems to possess considerable arbitrariness.

### [Phil] CREDO 1

I don't believe in the unity of mathematics and think that as time goes on the subject called mathematics becomes less and less unified. The Unity of Mathematics is a dream, a chimera, an ideal that doesn't exist:

- The *2010 Mathematics Subject Classification* (*MSC*) lists almost a hundred mathematical subjects. To some extent, each subject has its own techniques, intellectual resources and devotees. While there may, indeed, be some connections between e.g., potential theory and non-associative algebras and collaboration between experts that indicate a certain degree of unity and coherence in the field of mathematics I find the lack of unity more strikingly located elsewhere.
- *Diachronic and cross cultural disunity*. Written mathematics is easily 4000 years old. It has been created by people and has served for people a variety of purposes. A mathematician lives in a sub-culture at a certain time and place. A piece of mathematics does not exist only in a sequence of special symbols because the naked symbols are essentially uninterpretable. The symbols are embedded in a cloud of knowledge, meanings, associations, experiences, imaginations that derive from the *particularities* of time, place, person and the enveloping society.

---

[1] Contribution to the XI Österreichisches Symposion zur Geschichte der Mathematik, Organiser: Christa Binder, Topic: *Der Blick aufs Ganze. Gibt es große Linien in der Entwicklung der Mathematik*? Venue: Miesenbach, 22-28 April, 2012
[2] Email: <booss@ruc.dk>
[3] Email: <philip_davis@brown.edu>

- *Pythagoras* asserts that 3 is the first male number. In certain Christian theologies it is the number of the Godhead. If in Indian numerology the numbers 1,10,19, and 28 are *ruled by the sun*, the meaning of and the belief in those words may escape my readers. Historians of mathematics often explain a piece of ancient mathematics in terms of contemporary concepts. This may be anticipated because of the difficulty and ultimate impossibility, noted by numerous authors, particularly by ELEANOR ROBSON[4], of entering into the heads of the Past.[5]

- *Semantic ambiguity*. I may write down the sequence   x ╚∩σ∑≡ 6 and claim this is a piece of mathematics. But this claim cannot be substantiated on the basis alone of the mere symbols. To provide meaning, every mathematical statement must be embedded in a narrative in some natural language (English, German, et alii.) Furthermore, its significance as mathematics cannot be established if its knowledge is limited to one and only one person. (Private revelation.)

- *Semiotic ambiguity*. Can it be determined when two mathematical statements phrased differently, are asserting the same thing? BARRY MAZUR[6] has begun a discussion of this question.

- *Non-acceptance or doubts about certain theories put forth by professional mathematicians*. Examples are easily found. Originally there was one formal geometry: that of EUCLID. After BOLYAI and LOBACHEVSKII there were three; and after RIEMANN an infinity of geometries. ZERMELO did not believe GÖDEL's proof of the Incompleteness Theorem. For GEORGE BERKELEY: Infinitesimals were the *ghosts of departed quantities*. The skepticism of KRONECKER, POINCARÉ, ZERMELO, E. PICARD, BROUWER, HERMANN WEYL, WITTGENSTEIN, ERRETT BISHOPP ET ALII  regarding the concepts of CANTOR.

- A well known quote from the great applied mathematician RICHARD HAMMING sums it up:

---

[4] Eleanor Robson, *Mathematics in Ancient Iraq: A Social History*, Princeton and Oxford: Oxford University Press, 2008, xxvii + 472 pp., ISBN: 978-0-691-09182-2. On the other hand, Robson is fully aware of a close connection between the "algebra" her Babylonians talk about and what we can find in EUCLID's *Elements II* and in later algebra of equations.

[5] See André Weil, History of mathematics: why and how, Plenary Lecture, in: Proceedings of the International Congress of Mathematicians, Helsinki, 1978, pp. 227-236. He was most outspoken in his demand to explain a piece of ancient mathematics in terms of contemporary concepts, while at the same time warning against *anachronisms*: "An understanding in depth of the mathematics of any given period is hardly ever to be achieved without knowledge extending far beyond its ostensible subject-matter. More often than not, what makes it interesting is precisely the early occurrence of concepts and methods destined to emerge only later into the conscious mind of mathematicians; the historian's task is to disengage them and trace their influence or lack of influence on subsequent developments." Ironically, Weil's attitude was anticipated 130 years earlier by a German thinker: "The so-called historical presentation of development is founded, as a rule, on the fact that the latest form regards the previous ones as steps leading up to itself", Karl Marx, *Outline of the Critique of Political Economy* (*Grundrisse*), 1857-61, Penguin 1973, in elaborating his famous dictum *Human anatomy is a key to the anatomy of the ape*. See also the ontogenetic dilemma of unity: Nobody can think like an embryo. We can only describe our "adult" stage, even though the embryo is a preliminary stage of an "adult" stage.

[6] Barry Mazur, WHAT IS... a motive?, *Notices of the AMS* **51**/10 (2004), 1214-1216.

> I know that the great HILBERT said *We shall not be driven out of the paradise that CANTOR has created for us*, and I reply *I see no reason for walking in.*

- *Philosophic ambiguity*. Prior to the end of the 19th Century there was one philosophy of mathematics: that of Platonism. Now there are easily five distinguishable philosophies: together with variations that exhibit the Freudian *narcissism of slight differences*.

- *And yet....* There is something that is called mathematics. The history of the wool trade in 14th Century Brabant, cited by IBSEN, is not mathematics. *Sag' mir: Wo versteckt sich die Einheit der Mathematik?*

## [Phil] CREDO 2

- I believe that mathematics cannot have foundations and, in fact, doesn't need them.

## [Bernhelm] There is a general human longing for unity among great themes. More specifically there is a longing for them among mathematicians

I certainly see the points you make and to some extent I agree with you. Take, e.g., the *Duhem–Quine holism* thesis so popular in philosophy of science: The Thesis bids us to keep all things in view and argues strongly against the validation of single statements in isolation from all (!) possible connections. You and I have always agreed that, mildly speaking, these claims are utterly unrealistic. In mathematics and physics, the related discussions of the *Vienna Circle* have faded during the past 90 years.

However, as human beings, we need orientation and continuity. Regarding mathematics, perhaps CHARLES SANDERS PEIRCE was the philosopher who struggled most with the epistemological concepts of unity vs. independence. Basically, he pointed to the anthropogenic character of our thought concepts developed through hundred of thousands years, namely, our experience with procuring food and shelter and gaining a mate. These innate capacities should, however, be strengthened by a logically controlled abandonment of common sense views when confronted with phenomena beyond the shared human phylogenetic experiences. I believe that PEIRCE would like our *Blick aufs Ganze,* but also ÁGNES HELLER[7] and the late GIAN-CARLO ROTA[8] acknowledge a specific human habit when confronted with a variety of phenomena, namely a striving for expla-nation and meaning and the desire to memorize, communicate and reconsider the findings.

---

[7] Ágnes Heller, Can the unity of sciences be considered as the norm of sciences? In: Helga Novotny and Hilary Rose (eds.), *Counter-Movements in the Sciences – The Sociology of the Alternatives to Big Sciences*, D. Reidel Publ. Comp., Dordrecht, 1979, 57-66.

[8] Gian-Carlo Rota, *Indiscrete Thoughts*, Birkhäuser Boston, 1996.

This is why we have language, and though *human language* has throughout our history become both more specific and more diversified, nevertheless, there are great linguistic themes and, as against WHORF's Hypothesis that the structure of language affects the way we conceptualize things, there is a capacity for expressing a large variety of human observations and feelings in a shared way. Over the short time span of only 10,000 years, as archeologists tell us, dogs have devolved from their wolf ancestors into a remarkable variety of breeds; but in spite of all their current differences they still share characteristic features.

For me, the *power of mathematical formalisms*, can be derived from the capacity of formulae
- to recall and communicate condensed experiences and
- to suggest imaginative alternative approaches to already existing practices.

To exercise this cognitive transfer, the mathematical practitioner needs to discern and avail him/herself of the major themes in mathematics. Unity, the perception of unity and the search for unity, is constitutive for mathematics as a scientific subject.

## [Bernhelm]

Here is my view of the simultaneous tendencies of *specialization* (*diversification*) and *generalization* (*unification*) in mathematics:

1. Well-intended, but less well founded educational initiatives of the 1960s that centered elementary and advanced math teaching around sets and structures were readily overcome after 10 years of having been introduced. Math teachers on all levels rediscovered the challenges of teaching concrete mathematics. *Generalized Abstract Nonsense – GAN* -- was abolished.
2. Math research has shown a remarkable and powerful counter-movement against excessive generalizations. I mention here only four cases all of which are related to my own work in mathematics: (1) Returning to and the reconsideration of *generic cases* instead of striving for the greatest generality. (2) Orientation towards *algorithmic questions* under limited conditions instead of stating general non-viability. (3) Focus on *error quantities* such as the index, the eta-invariant, the spectral flow, the Maslov index. (4) Biology of *focused* systems, e.g., of a single cell instead of whole-body modeling holistically perceived.
3. Some university mathematicians experience such great pressure to publish, to plagiarize or blindly to resow in the same strip, that they feel they have not the time to think about the meaning of their work. Some are apt to inculcate the same snob feeling of high-standard accomplishment to their classes. This tendency is supported by a hierarchical division among mathematics in which a very few are the *architects*, full of seminal visions, and the many merely maintain the ground by filling in details or at best doing some *plumbing*.

4. In view of the experiences with the New Math, a threat to the intellectual unity of mathematics arises from *declamations* of unity that reduces mathematical dissemination to shallow definitions (e.g., vector spaces, groups, limits) hoping to create meaningful essentials but in the absence of meaningful applications.[9]

## [Phil]

A propos of *Major Themes,* a quotation from TOCQUEVILLE struck me as pertinent. Of course, TOCQUEVILLE was writing about systems of government and not about mathematics, but I think it might elicit a response from thoughtful mathematicians:

> Men of democratic centuries like general ideas because they exempt them from studying particular cases; they contain, if I can express myself so, many things in a small volume and give out a large product in a little time. When, therefore, after an inattentive and brief examination, they believe they perceive a common relation among certain objects, they do not push their research further, and without examining in detail how these various objects resemble each other or differ, they hasten to arrange them under the same formula in order to get past them.[10]

## [Bernhelm]

I like the moderate conservatism of the preceding quote. JACOB BURCKHARDT coined the phrase *terribles simplificateurs*. A German proverb of uncertain origin states *Der Teufel steckt im Detail – The devil is in the details.* So look carefully! Recently, a professor of musicology at Aarhus University pointed to a conflict between knowledge and theory, deploring that her students were much better at reading BOURDIEU than reading notes. BERTOLT BRECHT's Herr K. has this wonderful remark about the problem of clipping a laurel hedge into a ball: *Well, there is the ball now, but where is the laurel?* Indeed, there are good reasons to be alarmed whenever we are confronted with *Große Linien* and *Der Blick aufs Ganze*; and where does that leave us in our debate about the *Große Linien in der Entwicklung der Mathematik*?

---

[9] A model for such a well-intended and well-written, but somewhat misleading advocacy of *unifying and generalizing concepts* can be found in the widely read and cited Jean-Luc Dorier, Meta level in the teaching of unifying and generalizing concepts in mathematics, *Educational Studies in Mathematics* **29** (1995), 175-197. The author's favorite example of a unifying and generalizing concept is *Linear Algebra.* Based on an extensive historical study of a single aspect of the genesis of the abstract concept of a vector space (for J-L Dorier, it is only the concept of space and its algebraic formalization), he arrives at an epistemological analysis and an analysis of teaching sequences with many interesting observations, but where the core concept of linearity, as most mathematicians will see it, namely the concepts of eigenvalues and spectrum, are absent. For contrast, see Peter D. Lax, *Linear Algebra*, Wiley 1997, where he frankly admits the rather dullness of the axioms of linear algebra and then continues: "It is astonishing that on such slender foundations an elaborate structure can be built, with romanesque, gothic, and baroque aspects. It is even more astounding that linear algebra has not only the right theorems but the right language for many mathematical topics, including applications of mathematics." (l.c., p. 1).

[10] Alexis de Tocqueville, *De la démocratie en Amerique* (1835/1840)—*[Democracy in America](#)*. It was published in two volumes, the first in 1835, the second in 1840. Various English language versions. The French original and an English translation (by Henry Reeve) are on the web in public domain: http://fr.wikisource.org/wiki/De_la_d%C3%A9mocratie_en_Am%C3%A9rique and
http://ebooks.adelaide.edu.au/t/tocqueville/alexis/democracy/complete.html

## [Phil]

Here are some of the turning points in the history of mathematics that have had consequences in the philosophy of mathematics:

1. Pythagorean Theorem. Sqrt (2).  (Existence)
2.  EUCLID's *Elements*. (Axiomatics. Idealization)
3. Algebraization of Arithmetic circa 15th C.  (Formalization)
4.  Discovery of the complex numbers. (Semantics)
5. Algebraization of Geometry. DESCARTES. (Downgrading the visual)
6. Invention of Calculus. NEWTON, LEIBNIZ. (Existence of infinitesimals)
7. Algebra goes abstract. GALOIS, HAMILTON. (Formalization)
8. Mathematical logic. BOOLE, FREGE, RUSSELL, WHITEHEAD. (Logicism)
9. Non-Euclidean geometry. BOLYAI, LOBACHEVSKII. (Conflict between empiricism and axiomatics.)
10. Axiomatization of the real numbers and of analysis. CAUCHY, WEIERSTRASS, et al. (Formalization.)
11. Cantorian set theory. (Existence)
12. Space goes abstract. RIEMANN, KLEIN, PEANO, HILBERT. ( Formalism, Degradation of the visual)
13. HILBERT's Program. GÖDEL's Incompleteness Theorem. (Destruction of Logicism)
14. Electronic digital computing machines and the subsequent deep mathematizations of all aspects of society. Change in mathematical research methodologies. (Preeminence of the discrete over the continuous)
15. Increasing relevance of stochasticism. (Ontology)

## [Bernhelm]

In my view, a new type of unity emerged in the *Renaissance* with the dissolution of the sensus communis, the vanishing of the basically common language of επιστήμη and the subsequent universalization of the method of the *mathematization of the natural sciences*.

Just for a few seconds, let me play on the common pride of mathematicians regarding GALILEI's famous dictum:

> La filosofia é scritta in questo grandissimo libro che continuamente ci sta aperto innanzi a gli occhi (io dico l'universo), ma non si puó intendere se prima non s'impara a intender la lingua, e conoscer i caratteri, ne' quali é scritto. Egli é scritto in

> lingua matematica, e i caratteri son triangoli, cerchi, ed altre figure geometriche, senza i quali mezi é impossibile a intenderne umanamente parola; senza questi é un aggirarsi vanamente per un oscuro laberinto.[11]

In ÁGNES HELLER's characterization:

> So it was the new (symbolic) language of natural sciences that *became* the *sole sensus communis* in an age of dissolution of integrations, communities, other types of *sensus communis*, the sole scientific language whose norm it is that it could be spoken by every one and in an equal manner. This language has developed in an age in which the universal concept of humanity as abstracted from religion, race and nation was born.[12]

Today we rightly consider the *Glasperlenspiel* of mathematicians; the mathematization of the sciences and of technology; and a sober approach to international relations *more geometrico* (GROTIUS) for a triumph of humanity and not a regression.

*Insisting on human value and meaning.* Skeptical voices, however, appeared on the scene in parallel with the emphasis on unity. There is a considerable price, they argued, attached to this new common language, to this new conception of objectivity and science. In HELLER's words, it has to "pay the price of being abstracted from everything that is human, for the ever given societality, from value ideas of moral and non-moral type". She credits in particular KANT and HUSSERL for delineating the limits of natural sciences and emphasizes HUSSERL's thesis according to which "the emergence of modern natural sciences is an historical achievement; consequently their world-constitution is reversible".

Similarly, GIAN-CARLO ROTA, while recognizing the positive cultural and technological contributions of mathematics and mathematization fought a life-long battle against the *pernicious influence of mathematical thinking on philosophy* as exemplified by analytical philosophy, The outlines of the new meaning-oriented unity are not yet clearly drawn. The dominant philosophy of mathematics is still moving in the realm of GALILEO GALILEI's quote, as witnessed by the contributions to a conference on *The Unity of Mathematics*, held in 2003 in honor of I.M. GELFAND's 90th birthday.[13] The late GELFAND himself, however, called for greater awareness of ongoing changes of the content and role of mathematics and insisted on the distinction between meaningful and meaningless abstractions and constructions:

> We have a *perestroika* in our time. We have computers which can do everything. We are not obliged to be bound by two operations - addition and multiplication. We also have a lot of other tools. I am sure that in 10 to 15 years mathematics will be absolutely different from what it was before.[14]

---

[11] Galileo Galilei, *Il Saggiatore, Lettere, Sidereus Nuncius, Trattato di fortificazione*, in: *Opere*, a cura di Fernando Flora, Riccardo Ricciardi Editore, 1953. Here: Il Saggiatore, cap. 6. In English: ``Philosophy is written in that great book which ever lies before our eyes — I mean the universe — but we cannot understand it if we do not first learn the language and grasp the symbols, in which it is written. This book is written in the mathematical language, and the symbols are triangles, circles and other geometrical figures, without whose help it is impossible to comprehend a single word of it; without which one wanders in vain through a dark labyrinth." *The Assayer* (1623), as translated by Thomas Salusbury (1661), p. 178, as quoted in *The Metaphysical Foundations of Modern Science* (2003) by Edwin Arthur Burtt, p. 75.
[12] L.c., p. 59.
[13] P. Etingof, V. Retakh and I.M. Singer (eds.), *The Unity of Mathematics - In Honor of the Ninetieth Birthday of I.M.Gelfand*, Birkhäuser, Boston, 2006, XXII + 631 pages, ISBN-10 0-8176-4076-2, e-IBSN 0-8176-4467-9.

And

> An important side of mathematics is that it is an adequate language for different areas: physics, engineering, biology. Here, the most important word is adequate language. We have adequate and nonadequate languages. I can give you examples of adequate and nonadequate languages. For example, to use quantum mechanics in biology is not an adequate language, but to use mathematics in studying gene sequences is an adequate language.

The emergence of a new type of unity, oriented differently, may be sensed in the outspoken ethical stand of M.F. ATIYAH, another exponent of the classical GALILEAN quote. Firstly, as president of the Royal Society and later as president of the Pugwash nuclear disarmament movement, he blamed the development and the consequent degradation of much mathematics on its applications to war and to juke boxes.

## [Bernhelm] Promising offshoots and developments

There are clearly distinguishable mainstreams in mathematics. The active research mathematician has continuously to make a choice as to what are the prominent and promising fields to enter into or to rely on their own originality and inspiration. The difficult and often narrow problem of choice goes back to LAGRANGE who expressed very definitely his conviction that now all what could be solved in mathematics had been solved while at the same time opening wide fields of new mathematical research.[15]

In 1933, NORBERT WIENER characterized the *hierarchy* of mathematical objects:[16]

> In the hierarchy of branches of mathematics, certain points are recognizable where there is a definite transition from one level of abstraction to a higher level. The first level of mathematical abstraction leads us to the concept of the individual numbers, as indicated for example by the Arabic numerals, without as yet any undetermined symbol representing some unspecified number. This is the stage of elementary arithmetic; in algebra we use undetermined literal symbols, but consider only individual specified combinations of these symbols. The next stage is

---

[14] L.c., p.xiv.

[15] "There is but one universe, and it can happen to but one man in the world's history to be the interpreter of its laws." That is what Lagrange is quoted to have said about Newton, according to Th. Kuhn, The function of dogma in scientific research, in: A.C. Crombie (ed.), *Scientific Change*, Heinemann, New York, 1963, pp. 347-369, here p. 353. Kuhn's own comment: "In receiving a paradigm the scientific community commits itself, consciously or not, to the view that the fundamental problems there resolved have, in fact, been solved once and for all."

[16] Norbert Wiener, *The Fourier Integral and Certain of its Applications*, Cambridge University Press, Cambridge, 1933, p. 1.

> that of analysis, and its fundamental notion is that of the arbitrary dependence of one number on another or of several others -- the function. Still more sophisticated is that branch of mathematics in which the elementary concept is that of the transformation of one function into another, or, as it is also known, the operator.

Today, we might be inclined rather to make long lists of promising developments and major themes abandoned to illustrate the nature of contemporary productivity[17]. Since Lagrange's pronouncement of the victorious end of mathematics and its putative revitalization, there has been permanent and productive tension between what has been accomplished and what new theoretical insights might lead to new fields. Every time a question seemed to be settled and a new fact established, new concepts have arisen on a higher level of abstraction. BØRGE JESSEN once quoted to me a remark of HARALD BOHR that all developments require and receive consolidation: for example, invariants were consolidated in groups, equations in operator algebras, statistics in probability, optimization in functionals. Instead of the much feared atomization of mathematics, a world of cross connections has been discovered and elaborated. With hindsight it is incomprehensible why JOHN VON NEUMANN declined the invitation to the 1954 Amsterdam ICM to give a HILBERT style talk that would present a list of the most important and as yet unsolved mathematical problems. On the basis of his work for the US Atomic Energy Commission (ACE) he would have been the ideal witness for the ever and ever more manifest unity of mathematics. To me it seems that only regards to military security prevented him of demonstrating that it had become easier to oversee mathematics since HILBERT's 1900 and not more difficult and certainly not impossible, as VON NEUMANN claimed in his famous letter.

Underlying all specializations and generalizations, there is one dominant theme in the development of mathematics, namely, striving for meaning: for human meaning. Such meaning may be found in many directions, aesthetic, cognitive or utilitarian. To me, when all has been said, the search for, the discovery and the construction of meaning establish a kind of unity within mathematics.

## [Phil] The search for *Unity within Diversity* as a never ending process

There is certainly unity within mathematics. The Brown University catalog lists 50 different courses under one heading: Mathematics. Mathematicians of the world gather together every four years. On the other hand, Applied Math at Brown split off from Pure Math, and Computer Science split off from Applied Math.

I think that the phenomenon we are dealing with goes under the name of *Unity within Diversity*. This is a vast topic that spans all intellectual disciplines (Google the italicised phrase!) and the search for such unity within diversity is a never ending process.

---

[17] See Philip J. Davis, The rise, fall, and possible transfiguration of triangle geometry: A mini-history, *The American Mathematical Monthly* **102**/3 (March 1995), 204-214.